\documentstyle[10pt]{article}
\newcommand{\R}{{\bf R}}

\newcommand{\Z}{{\bf Z}}
\newcommand{\N}{{\bf N}}
\newcommand{\C}{{\bf C}}
\begin{document}

\title{Lazzeri's Jacobian of oriented compact riemannian manifolds}
\author{ Elena Rubei}
\date{ }
\maketitle

{\small {\bf E-mail address:} rubei@mail.dm.unipi.it}

{\small {\bf Math. classification:} Primary 14K30; Secondary  53C99.}

{\small {\bf Keywords:} Jacobians, riemannian manifolds.                 }

\vspace{0.4cm}

{\footnotesize {\bf Abstract}
The subject of this paper is  a Jacobian, introduced by F. Lazzeri
(unpublished),
associated to every compact oriented  riemannian manifold whose dimension 
is twice an odd number.
 We start the investigation 
of Torelli type problems and Schottky type problems for Lazzeri's Jacobian;
in particular we examine the case of tori with flat metrics.
Besides we study Lazzeri's Jacobian for K\"{a}hler manifolds and its 
relationship with other Jacobians. 
Finally we examine Lazzeri's Jacobian of a bundle.}

\section*{1. Introduction}

The subject of this paper is  a Jacobian, introduced by F. Lazzeri
(unpublished),
associated to every compact oriented riemannian manifold  whose
 dimension is twice an  odd number; if $(M,g)$
 is a compact oriented riemannian manifold  of 
 dimension $n=2m =2(2k+1)$,
Lazzeri's Jacobian of $(M, g)$ is the following p.p.a.v.:
                       
the torus $H^{m}(M,\R)/(H^{m}(M,\Z)/torsion)$ with the complex 
structure given by the operator $\ast$ and the polarization whose imaginary 
part is $(\alpha , \beta )=-\int \alpha  \wedge \beta $,  $\alpha ,\beta  \in 
H^{m}(M,\Z)$.

 (It is really a p.p.a.v., in fact:  the operator $\ast : H^{m}(M, \R) 
\rightarrow H^{m} (M, \R)$ 
(defined through the isomorphism of $H^{m}(M, {\bf R}) $ with the space of 
harmonic $m$-forms)
has square $-1$, so it induces a complex structure
on $ H^{m} (M, \R)$; besides  $ \int \cdot \wedge \cdot =
\int \ast \cdot \wedge \ast \cdot$ and  $ \int \cdot \wedge \ast \cdot $
is positive definite.)

One can easily see that 
the following definition is good:

{\bf Definitions 1.1.} {\it 
Let  $k \in {\bf N}$ and $m =2k+1$.
 Let ${\cal R} \subset \{(M,g) \, | \;  M$
  oriented compact $C^{\infty}$
manifold of dimension $2m$, $g$ riemannian 
metric on $M\} $. 
Let $\sim $ be the following equivalence relation on ${\cal R}$:
 $(M_{1},g_{1}) \sim (M_{2},g_{2})$ iff there are an orientation 
preserving diffeomorphism $f: M_{1} \rightarrow M_{2}$
and a $ C^{ \infty } $ map $t: M_{1} \rightarrow 
\R^{+} $  such that $(M_{1}, t \, g_{1}) \in {\cal R}$ and
 $f^{\ast}g_{2}=t \, g_{1}$. 
Let ${\cal A}_{h}$ be the set of 
p.p.a.v.'s of dimension $h$
 up to isomorphisms. 
Let 
\[ T_{\cal R}:
{\cal R}/\sim \;\; \longrightarrow \; {\cal A}_{1/2 \, b_{m}(M)} \]
be the map  associating to the class of $(M,g)$ 
Lazzeri'Jacobian of $(M,g)$.

Now  fix  $M$  
and  a symplectic basis  of $H^{m}(M,\Z)/torsion$ with respect 
$ - \int \cdot \wedge \cdot $ and          
let  ${\cal R} \subset \{g\;|\;g\;metric\;on\;M\}$. 
Then $ T_{{\cal R}}$ can be lifted   to a map                      
        \[ \hat{T}_{\cal R}:{\cal R}/ conformal\, equivalence
\;\longrightarrow \; {\cal H}_{ 1/2 \, b_{m}(M)}, \]                    
where ${\cal H}_{h}$  is the $h$-Siegel upper half space 
(two metrics on $M$, $g_{1}$ and $g_{2}$, are said conformal equivalent 
iff $g_{1} =t  g_{2}$ where $t$ is a 
$ C^{ \infty }$ map: $M \rightarrow \R^{+}$).

We  often omit the subscripts ${\cal R}$ in $T_{\cal R}$ and
$\hat{T}_{\cal R}$.}

If we consider the metrics of constant curvature equal to $-1$ on 
a $C^{\infty}$  compact surface, we have that Lazzeri's 
Jacobian is the   usual Jacobian of the correspondent Riemann surface. 

The outline of the paper is the following:
Section 2 deals with 
 Torelli and Schottky  problem for Lazzeri's Jacobian
 of tori with flat metrics,
in Sect. 3 we study Lazzeri's Jacobian for K\"{a}hler manifolds and its 
relationship with other Jacobians and 
in Sect. 4  we examine Lazzeri's Jacobian of a bundle.

\vspace{0.1cm}

 To state our results on Lazzeri's 
Jacobian we fix here some notation   we will use in all the 
paper.

{\bf Notation  1.2.} {\em
$\bullet$ {\boldmath ${\cal A}_{h}$, ${\cal H}_{h}$:}
defined above;
we will often drop the subscripts $h$. 

$\bullet$  {\boldmath $G_{\R}(h,V), G_{\C}(k,W)$.}
 $G_{\R}(h,V)$ is the Grassmannian of the real
 $h$-subspaces of the real vector space $V$,
 $G_{\C}(k,W)$ is the Grassmannian of the complex
 $k$-subspaces of the complex vector space $W$.

$\bullet$ {\boldmath${\cal F}_{n}$.}
We define ${\cal F}_{n}= \{(\R^{n}/\Lambda,g)\; | \; \Lambda \;\;
lattice, \;\; g \;\; flat \;\; metric \;\; on \;\;
 \R^{n}/\Lambda \}/\sim$,
where $\sim $ is the equivalence relation defined in Def. 1.1; here it 
becomes  
$(\R^{n}/\Lambda, g) \sim (\R^{n}/\Lambda', g')$ iff $\exists$
a  orientation preserving map
 $f: \R^{n}/\Lambda \rightarrow 
\R^{n}/\Lambda'$ induced by a linear map $ \R^{n} \rightarrow \R^{n} $ 
such that $f^{\ast}g'=c g$ for some $c \in \R^{+}$, where 
the orientation is the standard one of $\R^{n} $, (in fact        
  if $(\R^{n}/ \Lambda, g) $ and $
(\R^{n}/ \Lambda', g')$ are equivalent for $\sim $  through 
a  map $\varphi$, then $\varphi $ is given by an affine map $\R^{n} 
\rightarrow \R^{n}$).

$\bullet$ {\boldmath${\cal T}_{n}$.} 
Set    ${\cal T}_{n} := \{ T \in M(n \times n,  \R) \mid  T  \; lower 
\; triangular, \; T_{i,i} >0 \; \forall i , \; det \, T = 1  \} $.

$ \bullet $ Let $m \in \N$ and $n=2m$.
Consider the standard basis of $\R^{n}$, $\{e_{i}\}_{i=1,..n}$.
We call {\bf lexicographic ordered basis} of  
$\wedge^{m} \R^{n} $ the following ordered basis:
$e_{I}$, $ I= (i_{1},..., i_{m})  \in {\N}^{m}$  with $ 1 \leq 
i_{1} < i_{2} < ... < i_{m} \leq n $, 
  with the multindeces ordered by the  lexicographic order.
Let ${\cal I} : = \{ I = (1,  i_{2} ,.., i_{m} ) \in  \N^{m}
 \,|\; 1< i_{2} < ..< i_{m}
 \leq n \}$. If $I \in {\cal I}$, we  choose and fix forever 
one of the  multindeces $J$ such that we obtain 
$ (I, J)$ from $(1,....,n)$ with an even number of
transpositions and we 
 call it  $\tilde{I}$. Let $ {\cal E} = \{ \tilde{I}\, |\, 
I \in {\cal I} \}$.
We want to define another ordered basis of  
$\wedge^{m} \R^{n} $: 
  we take first the multindices  in ${\cal I}$  in the lexicographic 
order: we call them  $I_{1}, I_{2},...$;
 then we consider  $\tilde{I}_{1}, \tilde{I}_{2},...$;
  we call {\bf the symmetric lexicographic ordered basis}
of $\wedge^{m} \R^{n} $ the ordered basis  
$e_{I_{1} } , e_{I_{2}}, ...,e_{\tilde{I}_{1} } , e_{\tilde{I}_{2}}, ...,$
and we call ``the symmetric lexicographic order'' the order 
 $I_{1}, I_{2},...,\tilde{I}_{1}, \tilde{I}_{2},...$ of the multindeces
of ${\cal I} \cup {\cal E}$. Analogously for $({\bf R}^{n}) ^{\vee}$ instead
of $\R^{n}$.

$\bullet$ Let $M$ be a complex manifold.
 If  $H^{q}(M,\C)= V \oplus W $, $\pi_{V,W}$  will
 denote  the projection onto $V$.
We will often omit the subscript $V,W$ in $\pi_{V,W}$ when it will be clear.

$\bullet$ {\boldmath 
$ K(M,g)$, $K'(M,g)$, $ F^{p}_{q}(M)$, $J_{q}(M)$.}                    
Let  $(M,g)$ be a compact K\"{a}hler
manifold of complex dimension $m$; let $\Omega $ be the 
(1,1)-form associated to $g$.

{\it $K= K(M,g):=\{[\eta] \in H^{m}(M,\C)| \,
[ \eta] = \sum_{r \geq 0,\, r  \equiv -m(m+1)/2 \; (mod \,2)}
[\Omega^{r} \eta_{r}] $ with $\eta_{r}$ primitive form  of degree 
$m-2r  \}$ and

$K'= K'(M,g):=\{[\eta] 
\in H^{m}(M,\C)| \, [\eta] = \sum_{r \geq 0, \,r 
\equiv 1-m(m+1)/2 \; (mod \,2)} [\Omega^{r} \eta_{r}] $
 with $\eta_{r}$ primitive form of degree $m-2r \}$}.

$ F^{p}_{q}(M)  :=  \oplus_{a+b=q,\; a \geq p} H^{a,b}(M,\C)$      
(we will omit the subscript $q$ when no confusion can arise).

    If $q$ is an odd number, we define
$ J_{q}(M)$ $ = \oplus_{a+b=q,\; a-b \equiv 1 \; (4)} H^{a,b}(M,\C)$.}

\vspace{0.2cm} 
                               
First we examine the case of flat metrics on tori.
  Let $n=2m= 2(2k+1)$ and $ N := \frac{1}{2} \dim \wedge^{m}\R^{n}$; 
we  study  $T:{\cal F}_{n} \rightarrow {\cal A}_{N}$ 
and, choosing a symplectic basis, we study $\hat{T}$ and 
in  particular Im $\hat{T}$ (and then Im $T$):

{\bf Theorem A. }                    
{\it {\bf i)} Choosen any symplectic basis of $H^{m}(\R^{n}/\Z^{n},\Z)$
with respect to $ -   \int \cdot \wedge \cdot$, the map 
$\hat{T}:   \{g\;|\;g\;flat\; metric\;on\;\R^{n}/\Z^{n}\}
/ conf.\; equiv . \;\longrightarrow \;{\cal H}_{N} $
is injective.
     
{\bf ii)} Now
call $\{ dx_{1},..., dx_{n}\}$
 the standard basis of $(\R^{n})^{\vee}$.
If we take  as a symplectic basis of $H^{m}(\R^{n}/\Z^{n},\Z)$
with respect to $ -   \int \cdot \wedge \cdot$, 
the sym. lex. ordered basis 
$\{ dx_{I} = dx_{i_{1}} \wedge ... \wedge dx_{i_{m}} 
\}_{I \in {\cal I} \cup {\cal E}}$, we have that  
$Im  \hat{T}= 
\{ X + i Y \in {\cal H}_{N} \mid 1) \; and \; 2) \; hold \}$, where:

1) \(Y= E E^{t}\) where 
 $E_{I,J}= det (T)_{I,J} $, $I,J \in {\cal I}$, 
${\cal I}$  ordered in  lex. order., for some $ T $ 
upper triangular matrix $ n \times n $
 with determinant $1$ and $ T_{i,i} \geq 0$    
(in  particular $E$ is upper triangular 
with positive diagonal elements and the entries of every column 
of $E$ are the Pl\"{u}cker coordinates of an element of 
$G_{{\bf R}}(2k, {\bf R}^{n-1})$ and the same for the entries of every
 row of $E$).

2) $  X_{IJ}= 0$  if  $I$ and $ J$  have more than 
one index in common,
where  $I,J \in {\cal I} $, ${\cal I}$  
ordered in  lex. order.         

{\bf iii)}   Lazzeri's Jacobian of a generic flat oriented torus has not
nontrivial automorphisms as p.p.a.v..
           
{\bf iv)} The map $T: {\cal F}_{n} \longrightarrow
{\cal A}_{N}$ is generically locally injective, {\bf v)} but not injective. }

In Section 3 we consider Lazzeri's Jacobian for K\"{a}hler manifolds $(M,g)$
and we
see that it depends only on the complex structure of $M$  and on the
cohomology class of the $(1,1)$-form associated to $g$; using  that, we 
examine  the relationship of Lazzeri's Jacobian with Weil's  and 
Griffiths' Jacobians:

{\bf Theorem  B.} {\it
Let $(M,g)$ be a compact K\"{a}hler manifold
of complex dimension $m=2k+1$; (in the sequel
$\pi $  is $ \pi_{J_{m}(M), \overline{J_{m}(M)}}$);
we have that, as complex tori,
                          
$\bullet $ $k^{th}$ Weil's Jacobian is 
$H^{m}(M,\R)/(H^{m}(M,\Z)/torsion)$ with 
the complex structure given by $C$
(or $J_{m}(M)/ \pi(H^{m}(M,\Z))$ with the complex structure given by 
$i$),      
                                  
$\bullet $ $k^{th}$ Griffiths' Jacobian is
 $H^{m}(M,\R)/(H^{m}(M,\Z)/torsion)$ with 
the complex structure given by $C$ on the set $\{\eta + \overline{\eta }\mid\;
\eta \in J_{m}(M) \cap \overline{F^{k+1}(M)}\}$ 
and $-C $ on $\{\eta + \overline{\eta } \mid \;
\eta \in J_{m}(M) \cap F^{k+1}(M)\}$
(or $J_{m}(M)/ \pi(H^{m}(M,\Z))$ with the complex structure 
$i|_{J_{m}(M) \cap   \overline{F^{k+1}(M)}}
 \oplus -i|_{J_{m}(M)\cap   F^{k+1}(M)}$), 
                  
$\bullet $
Lazzeri's Jacobian is $H^{m}(M,\R)/(H^{m}(M,\Z)/torsion)$ with 
the complex structure  $C|_{H^{m}(M,\R) \cap K} \oplus 
-C|_{H^{m}(M,\R) \cap K'}$
(or $J_{m}(M)/ \pi(H^{m}(M,\Z))$ 
with the complex structure  $i|_{J_{m}(M) \cap K}
\oplus -i|_{J_{m}(M) \cap K'}$);             
     
$k^{th}$ Weil's Jacobian and Lazzeri's Jacobian
are p.p.a.v.'s (also if the 1-1 form associated to the metric isn't rational)
and  the real part of the polarization is the same (equal to $\int \cdot
\wedge \ast \cdot$ on $H^{m}(M,\R)/(H^{m}(M,\Z)/torsion)$).}

Besides                                          
we consider another class of Ricci-flat metrics (besides the flat ones on tori)
 K\"{a}hler-Einstein metrics on the complex manifolds with trivial 
canonical bundle and  we find another  local Torelli 
theorem  for Lazzeri's Jacobian ({\bf Corollary C}),  and
one could formulate the following conjecture:

{\bf Conjecture.} {\it $\hat{T}_{\cal R}$ is locally
 injective if ${\cal R} $ is a set of Ricci-flat metrics on a manifold $M$.} 
                                
Finally in Section 4 
 given a bundle $F  \rightarrow M$, we study the relationship 
between Lazzeri's Jacobian of $F$ and Lazzeri's Jacobian of $M$
{\bf (Prop. D)}. 
              
We think 
that one can easily find  open problems about Lazzeri's Jacobian,
for instance to go on with the study of Schottky and Torelli type problems,
 to study Prym-Tyurin varieties for Lazzeri's Jacobians and  
the relationship between Lazzeri's Jacobian and the theory of degeneration of 
abelian varieties or more precisely 
to study  a possible ``object'' $T(M,g_{0})$ to 
associate to every compact oriented manifold $M$  of dimension 
$2(2k+1)$ with a singular metric
$g_{0}$ (i.e.   $(g_{0})_{P}$ is 
semipositive definite $ \forall \, P \in M$), 
such that, if $g_{t}$  are riemannian 
metrics on $M$ for $t>0$
 and $g_{t} \longrightarrow g_{0}$ for $t \longrightarrow 0$,
 then $T(M,g_{t}) \longrightarrow T(M,g_{0})$ in some sense.
Observe that if we consider  a torus $M= \R^{2(2k+1)}/ \Lambda $
 with a singular flat metric $g_{0}$, one could define $T(M,g_{0})$ in the
 following way:             
let $g$ be a flat metric  on $M$ extending  $g_{0}$ (i.e. 
  $(g_{0})_{P}= (g)_{P} $
on $ (\ker (g_{0})_{P})^{\perp_{g}}$ $\forall P \in M$, 
or equivalently for one point of $M$);  
we define  $T(M, g_{0})= {\cal H}^{2k+1}_{g}(M, \R)/ 
\psi ({\cal H}^{2k+1}_{g}(M, \Z))$ with the complex structure $\ast_{g}$,
where ${\cal H}^{2k+1}_{g}$ is the set of
  harmonic $(2k+1)$-forms for $g$,
i.e. invariant $(2k+1)$-forms,  
$\psi $ is defined pointwise as the projection  onto       
$ \wedge^{2k+1}( (\ker (g_{0})_{P})^{\perp_{g}})^{\vee}$ $ \forall  P$ 
and $\ast_{g}$ the  operator $\ast $ for $g$ 
(one can easily see that this definition doesn't depend on the metric $g$ 
extending $g_{0}$, see Appendix.)

\section*{2. The case of tori with flat metrics}

\subsection*{2.a. Some lemmas of linear algebra}

{\bf Lemma 2.1.} {\it 
 Let $ \Omega $ be a real upper  triangular matrix 
 $n \times n$, where $n=2m$. In the lex. ordered  basis of 
$\wedge^{m} \R^{n} $, the matrix   
$\bigwedge ^{m} \Omega $ is upper triangular. 
Besides, if  $ \bigwedge ^{m} \Omega $ is equal to    
{\footnotesize\( \left( \begin{array}{cc}
                                                  A & C \\
                                                  B & D         
\end{array} \right) \)}
in the sym. lex. ordered  basis of $\wedge^{m} \R^{n} $,
then     $B=0$, $A$ is  upper triangular and $D$ lower triangular.}

{\sc Proof.} Left to the reader. \hfill $\Box$

{\bf Lemma 2.2.} {\it  Let $ \Omega $ be a real upper triangular matrix 
 $n \times n$, where $n=2m$. 
       
Let us fix an ordered  basis of $\wedge^{m} \R^{n} $
defined in the following way:
order in any fixed way the multindeces
 $I  $ in  $  {\cal I}$; call them $K_{1}, K_{2},..  $,
 then, order the multindices in $ {\cal E} $ in the following way:
$\tilde{K}_{1}, \tilde{K}_{2},.. $; consider 
$e_{K_{1}} , e_{K_{2} },...,e_{\tilde{K}_{1}} , e_{\tilde{K}_{2}},...$; 
(for instance the sym.  lex. ordered basis). 
                        
In this ordered basis of $\wedge^{m} \R^{n} $
 let $ \bigwedge ^{m} \Omega $ be
{\footnotesize\( \left( \begin{array}{cc}
                                                  A & C \\
                                                  B & D         
\end{array} \right). \)}               
We have: 
          
a)$A^{t}D= (\det\,\Omega ) I$;
               
b)suppose $\det \Omega \neq 0$; for $I \in {\cal I}$ and $J \in  {\cal E}$}
             
\[ (A^{-1}C)_{IJ}= \left  \{ \begin{array}{ll}
0                        & if \; \tilde{I}\; and\; J\; have\; more\; than 
\;one\;index\; in\; common \\
\pm \frac{\Omega_{1,l}}{\Omega_{1,1}} 
& if\; \tilde{I} \;and\; J \;have\; only\; l\; 
in\; common
\end{array}     
    \right    . \]

{\sc Proof.}    
If $P$ and $Q$ are two multindeces, $\Omega_{P,Q} $ will denote 
                   the determinant of the minor $P,Q$ of $\Omega$.
We will denote  the $j^{th}$-column of $\Omega$ by $v_{j}$ and 
$v_{j}-\Omega_{1,j}e_{1}$ by  $\overline{v}_{j}$.

{\bf  a)}   
     Let   $I=(i_{1},...,i_{m}) \in {\cal I}$
and  $J=(j_{1},...,j_{m}) \in {\cal E}$. 
 Observe that

$(A^{t}D)_{I,J} e_{1} \wedge ...\wedge e_{n}
=\sum_{S \in {\cal I} }
  \Omega_{S,I} \Omega_{\tilde{S},J} e_{1} \wedge ...\wedge e_{n}=$

$= ( \sum_{S \in {\cal I} }
  \Omega_{S,I}  e_{s_{1}} \wedge ....\wedge e_{s_{m}}) \wedge
(\sum_{T \in  {\cal E}}
 \Omega_{T,J} e_{t_{1}} \wedge ...\wedge e_{t_{m}})= $

$= ( \sum_{S \in {\cal I} \cup {\cal E}}
  \Omega_{S,I}  e_{s_{1}} \wedge ....\wedge e_{s_{m}}) \wedge
(\sum_{T \in {\cal I} \cup {\cal E}}
 \Omega_{T,J} e_{t_{1}} \wedge ...\wedge e_{t_{m}})=$

$=     v_{i_{1}} \wedge ...\wedge v_{i_{m}} \wedge v_{j_{1}} \wedge ...\wedge 
v_{j_{m}}   =
\left\{ \begin{array}{ll}
     0  & \mbox{ if $J \neq \tilde{I}$ } \\
 (\det\, \Omega) \, e_{1} \wedge ...\wedge e_{n} & \mbox{ if $J = \tilde{I}$ }
       \end{array}
\right . $

where $S=(s_{1},...,s_{m})$, $T=(t_{1},...,t_{m})$,
and the last but one  equality holds because $ \Omega_{S,I}= 0 $ for 
$ S \in {\cal E}$ and $I \in {\cal I}$.

{\bf b)}                     
Let    $I, J \in {\cal E}$.
We remark that, since neither $I$ nor $J$  contain $1$, $I$ and $J$ have some
 index in common.    Let  $I=(i_{1},...,i_{m})$ and
 $J=(j_{1},...,j_{m})$  with $l =i_{r}=j_{s}$.
We have:
\begin{eqnarray*}
& & 
(\det \Omega) (A^{-1}C)_{IJ} e_{1} \wedge ... \wedge e_{n} \stackrel{a)}{=}
 (D^{t}C)_{IJ} e_{1} \wedge ... \wedge e_{n}=
 (\sum_{K  \in { \cal I} }  \Omega _{\tilde{K},I} 
\Omega _{ K,J})e_{1} \wedge ... \wedge e_{n} =  \\ 
& = &\overline{v}_{i_{1}} \wedge ...\wedge \overline{v}_{i_{m}} \wedge (
v_{j_{1}} \wedge... \wedge v_{j_{m}}- \overline{v}_{j_{1}} \wedge... \wedge
\overline{v}_{j_{m}}) 
 =  \overline{v}_{i_{1}} \wedge ...\wedge \overline{v}_{i_{m}} \wedge 
v_{j_{1}} \wedge... \wedge v_{j_{m}}  =  \\
& = & \overline{v}_{i_{1}} \wedge ...\wedge \overline{v}_{i_{m}} \wedge 
v_{j_{1}} \wedge...\wedge v_{j_{s-1}} \wedge  \Omega_{1,l} e_{1}
\wedge v_{j_{s+1}} \wedge ... \wedge v_{j_{m}}  = \\
& = &  \overline{v}_{i_{1}} \wedge ...\wedge \overline{v}_{i_{m}} \wedge 
\overline{v}_{j_{1}} \wedge...\wedge \overline{v}_{j_{s-1}} \wedge  
\Omega_{1,l}e_{1}
\wedge\overline{v}_{j_{s+1}} \wedge ... \wedge \overline{v}_{j_{m}} =  \\
& = &  \left\{ \begin{array}{ll}
0     & if \, I\, and\, J\, have\, more\, than \,one\,index \, in\, common \\
\epsilon \, \frac{\Omega_{1,l}}{\Omega_{1,1}} \, \det \, \Omega \,
 e_{1} \wedge ... \wedge e_{n}                  & 
if\, I \,and\, J \,have\, only\, l\, 
in\, common
\end{array}          \right .
\end{eqnarray*}          
where $\epsilon $ is the sign of the permutation taking $(i_{1},...,i_{m},
j_{1},...,j_{s-1},1,j_{s+1},...,j_{m})$ in $(1,...,n)$.
\hfill $ \Box $

{\bf Lemma 2.3.} {\it
Let $m \in {\bf N} $ and $n=2m$.
Consider  the set  of positive scalar products on $R^{n}$ up to conformal 
equivalence. 
The map, defined on this set, associating to the class
of a positive definite scalar 
product its operator  $\ast $ on $\wedge^{m} (\R^{n})^{ \vee }$, 
is injective.}

{\sc Proof.}                             
 We observe that,
 given a scalar product on  $\R^{n}$, two elements of $\R^{n}$, $v$ and
 $w$, are perpendicular iff there exists a  $m$-subspace of $\R^{n}$
perpendicular to $w$ and containing $v$; thus $v$ and
 $w$ are perpendicular iff $\exists \; \alpha \in \wedge^{m-1} \R^{n}$, 
$\alpha $ simple (i.e. of the kind: $v_{1} \wedge ...\wedge v_{m-1}$)
such that $w \wedge \ast (\alpha \wedge v)=0 $  and $ \alpha \wedge v \neq 0$.
Thus $\ast$ determines the conformal structure.
 \hfill $\Box$

\subsection*{2.b. Proof of Theorem A}

{\bf Remark 2.4.} {\it 
The set    $\{(\R^{n}/ \Z^{n},g) \; |  \; g \; 
  flat \;  metric \; on \; \R^{n}/\Z^{n}\}/ conf. \, equivalence\;$
 is in bijection with the set  $P_{n}$ 
of symmetric positive definite matrices 
$n \times n$ with determinant $1$,  thus with $ {\cal T}_{n}$.
         The set    ${\cal F}_{n}= \{(\R^{n}/ \Z^{n},g) \; |  \; g \; 
  flat \;  metric \; on \; \R^{n}/\Z^{n}\}/ \sim \; = 
P_{n}/ SL(n, {\bf Z})$ (where $A \in SL(n, {\bf Z})$ acts on $P_{n}$  by
$P \mapsto A^{t} P A $) and we endow it with the quotient topology induced by 
the set above.}

\vspace{0.1cm}

{\bf Proposition 2.5.} {\it
Let  $n := 2m= 2(2k+1)$ 
and $N:= \frac{1}{2} {\footnotesize \left( \begin{array}{c}
                                                 n  \\
                                                 m  
\end{array} \right)}$. 
      
Let $ L \in GL(n, {\bf R})$;        
let $\wedge ^{2k+1} (L^{-1})^{t}  =$
{\footnotesize $ \left( \begin{array}{cc}
                                                  A( L) & C( L) \\
                                                  B( L) & D( L)  
\end{array} \right)$}
in the  sym. lex. order;
set
{\footnotesize $ \left( \begin{array}{c}
                                               X( L)  \\
                                               Y( L)       
\end{array} \right) $} $=$
{\footnotesize  $ \left( \begin{array}{cc}
                                              A( L) &  -B( L) \\
                                              B( L) &   A ( L)      
\end{array} \right)^{-1} $}
{\footnotesize $ \left( \begin{array}{c}
                                               C( L) \\
                                               D( L)       
\end{array} \right); $}
define $Z( L)=X( L)+iY( L) \in M(N \times N, \C )$.

  Let   $det \, L > 0$.             
Call $\{ dx_{1},..., dx_{n}\}$
 the standard basis of $(\R^{n})^{\vee}$.
Choose the sym. lex. ordered basis 
$\{ dx_{I} = dx_{i_{1}} \wedge ... \wedge dx_{i_{2k+1}} 
\}_{I \in {\cal I} \cup {\cal E}}$
 as symplectic basis of 
$H^{m}({\bf R}^{n}/ {\bf Z}^{n},{\bf Z})$ for 
$ - \int \cdot \wedge \cdot$. 
We have that $Z(L) =  \hat{T} (\R^{n} /  \Z^{n}, L^{t}L)$
(the orientation of $ \R^{n} / \Z^{n} $
is given by the standard one of $\R^{n}$).}

{\sc Proof.}      
Let $b_{i}$ be  the columns of  $L^{-1}$, they are an orthonormal basis
of $ {\bf R}^{n}$ for the metric $ L^{t} L$.

The matrix representing the ordered basis
$\{dx_{I}\}_{I \in {\cal I} \cup {\cal E}}$  in function of 
the sym. lex. ordered basis of 
$H^{m}({\bf R}^{n}/ {\bf Z}^{n},{\bf R})$,
$\{b_{I}^{\vee}\}_{I \in {\cal I} \cup {\cal E}}$
 is $\wedge ^{2k+1} (L^{-1})^{t}=$
{\footnotesize $ \left( \begin{array}{cc}
                                                  A & C \\
                                                  B & D        
\end{array} \right)$}.                 
Then the matrix expressing   
$\{dx_{I}, \ast dx_{I} \}_{I \in {\cal I}}$  
   in function of 
the sym. lex. ordered basis 
$\{b_{I}^{\vee}\}_{I \in {\cal I} \cup {\cal E}}$
is 
{\footnotesize $ \left( \begin{array}{cc}
                                              A  &  -B \\
                                              B  &   A       
\end{array} \right)$}.  
Thus the matrix expressing $\{dx\}_{I \in  {\cal E}}$ 
  in function of  $\{dx_{I}, \ast dx_{I} \}_{I \in {\cal I}}$  
is:
{\footnotesize $ \left( \begin{array}{cc}
                                              A  &  -B \\
                                              B  &   A       
\end{array} \right)^{-1} $}
{\footnotesize $ \left( \begin{array}{c}
                                               C \\
                                               D       
\end{array} \right). $}                                      
   \hfill $\Box$

\vspace{0.1cm}

{\sc Proof of Theorem A.} In proving ii), iii) and v) we use Prop. 2.5
 and its notation.

{\bf i)}           
The map $\hat{T}$ is injective  by  Lemma 2.3 and because
the harmonic forms on tori are the translations-invariant forms.

{\bf ii) }
Observe that, if $L \in {\cal T}_{n}$, then
 $ \wedge ^{2k+1} (L^{-1})^{t} =
  {\footnotesize \left( \begin{array}{cc}
                            A(L) &  C(L) \\
                               0 &   (A(L)^{t})^{-1}
\end{array} \right)}$
by Lemmas 2.1 and 2.2. 
              
Then
 $X(L)= A(L)^{-1} C(L)$ and 
 Lemma 2.2 implies the  claim for $X$ .                        
               
Besides  $Y(L)=  A^{-1}(L) (A^{-1}(L))^{t}$,
where  $A(L)_{I,J}= \det ((L^{-1})^{t})_{I,J} $;
 then we have
$Y(L)=  A^{-1}(L) (A^{-1}(L))^{t} =  A(L^{-1}) (A(L^{-1}))^{t}$
where $A(L^{-1})_{I,J}= \det (L^{t})_{I,J} $.
Taking $ E = A(L^{-1}) $ and $ T = L^{t} $ we conclude.

{\bf iii)} Let $G$ be the modular group.
If $F_{\sigma} := \{ x \in {\cal H}_{N} \mid \sigma (x) =x \} $
 for   $\sigma \in G^{\ast } := G - \{ Id \}$,
we have to prove that ${\cal T}_{n}- Z^{-1}(\cup_{\sigma \in G^{\ast }}
 F_{\sigma})$ is a open dense subset of ${\cal T}_{n}$.
The openness follows from the fact  that 
$G$  acts properly and discontinously on  ${\cal H}$
(see \cite{L-B} p. 218).
To prove the density, it is sufficient (by Baire's theorem) 
to  prove that $ \forall \sigma  \in G^{\ast}$, 
the set ${\cal T}_{n} - Z^{-1}( F_{\sigma})$ 
is a open dense subset of ${\cal  T}_{n}$; since
$Z^{-1}( F_{\sigma})$ is defined by
polynomial equations,  we have only to prove that 
${\cal T}_{n} - Z^{-1}( F_{\sigma}) \neq \emptyset$ i.e.
that $Z({\cal T}_{n}) \subset F_{\sigma}$ implies $\sigma = Identity$.

Let $\sigma \in G $ be the  map $Z \mapsto (MZ+N)(PZ+Q)^{-1}$.
Obviously $Z \in F_{\sigma}$  iff 
 $ZPZ+ZQ=MZ+N$.

Let $L \in {\cal T}_{n}$; for every $c \in \R^{+}$, let 
$\tilde{L}(c)$ be the matrix obtained from $L$
by  multiplying $L_{1,1}$
by $c^{-4k-1}$ and the other entries of $L$ by $c$; 
 we have  $Z(\tilde{L}(c))=c^{4k+2}Z(L)$; since
  $Z(\tilde{L}(c) ) \in F_{\sigma}$ $\forall c \in \R^{+}$,
 we   have  
$c^{8k+4} Z(L)PZ(L)+c^{4k+2}(Z(L)Q-MZ(L))-N=0$ 
 $\forall c \in \R^{+} $, $ \forall L \in {\cal T}_{n}$;
 thus: 

1) $N=0$  \hspace{1cm}   
2) $Z(L)PZ(L)=0$  $\forall \;  L \in {\cal T}_{n}$ \hspace{1cm}
3) $Z(L)Q-MZ(L)=0$ $\forall \; L \in {\cal T}_{n}$ 
  
Hence:  $Q=(M^{t})^{-1}$ by 1) and  $\sigma \in G $;
 $P=0$, $M=Q$ by 2) and 3), taking $L=I$. Thus $M$ is 
orthogonal.
From 3), we have: $Y(L)M=MY(L)$ and  this  
implies $M$ diagonal (take  $L $ in such way that
$Y(L)$ diagonal with the elements
 of the diagonal different one from another). 
Being $M$  orthogonal and diagonal,
  $M$ must be  diagonal  with only $\pm 1$ on the diagonal;  
still from $Y(L)=MY(L)M^{-1}$, taking $L$ in such way that
$Y(L)$ not diagonal, we obtain $M=I$.
          
{\bf iv)} It follows right away from i) and iii), since       
  $G$ acts properly and discontinously on  ${\cal H}$.

{\bf v) } Let us consider two diagonal matrices
$\in {\cal T}_{n}$, one inverse of the other:
$F$ and $F^{-1}$.
        
We have: $ C(F)   =  B(F) = C(F^{-1})   =  B(F^{-1}) =0  $  and 
$ A(F)   =  A(F^{-1})^{-1} $,  $D(F)=  D(F^{-1})^{-1}$. 
Thus $Y(F^{-1})=Y(F)^{-1}$ and $X(F^{-1})=X(F)=0$.
Thus $Z(F)=-Z(F^{-1})^{-1}$; so the $Z$'s differ by a modular map.
Thus Lazzeri's Jacobians of $ (\R^{n} /  \Z^{n},  F^{2})$ and
of $ (\R^{n} / \Z^{n} , (F^{-1})^{2})$ are isomorphic p.p.a.v.'s.
But  we can choose  $F$ in such way  that 
$\not \exists  A \in SL(n, {\bf Z}) $ s.t. $ A^{t} F^{2} A = (F^{-1})^{2} $.
 \hfill $ \Box $

\section*{3. The case of K\"{a}hler manifolds}

   We recall some facts on the  operators $ \ast $, $ C $ and $L$ 
(see for instance  \cite{Ch}, \cite{Wei}, \cite{G-H}, \cite{Wel}, \cite{Gre}).

\vspace{0.1cm}

Let  $M$ be a complex manifold of complex dimension 
$m$; Weil's operator $C$ is defined on 
the forms of bedegree $(a,b)$ by $C \eta = i^{a-b} \eta $ and then  is 
extended by linearity to arbitrary forms;  observe that $C$ 
 takes real forms to  real forms.
If $M$ is  compact 
and of  K\"{a}hlerian type,  then, $ \forall q \in \N$ odd,
the operator $C$ defines a complex structure on $H^{q}(M,\C)$ (and 
thus on $H^{q}(M,\R)$) depending only on the complex structure of $M$; 
the $i$-eigenspace of $C$ is 
$J_{q} (M) = \oplus_{a+b=q, a-b \equiv 1 \; (4)} H^{a,b}(M,\C)$.

\vspace{0.1cm}

Now suppose  that $(M,g)$ 
 is a hermitian manifold and let
$\Omega $ be the $(1,1)$-real form associated to $g$; let
 $L$ be the operator on the set of forms defined by 
 $L \eta = \Omega  \wedge \eta $ and $ \Lambda $ 
the adjoint operator of $L$. 
A form $ \eta $ such that $ \Lambda  \eta =0 $ is said primitive.
Every form $\omega $ of degree $q$ can be written uniquely in the form
$\omega = \sum_{r \geq (q-m)^{+}} L^{r} \omega_{r}$, where $\omega_{r}$ is a
primitive form of degree $q-2r$.      
If  $g$ is of K\"{a}hler
then $L$ defines an operator on $  H^{q}(M, \C)$; 
 a class $[\omega] \in H^{q}(M,\C)$
is said primitive if $[\Lambda \omega]= 0$ and a class $[\omega ] \in 
H^{q}(M,\C)$ can be written uniquely in the form
$[\omega] = \sum_{r \geq (q-m)^{+}} L^{r} [\omega_{r}]$, where $[\omega_{r}]$
 is a  primitive class of degree $q-2r$.

\vspace{0.1cm}
               
{\bf Warning 3.1.} Here the operator $\ast 
$, defined as usual  on  real forms, is extended to complex forms by 
$\C$-linearity as in \cite{Ch} and \cite{Wei}.

\vspace{0.2cm}

{\bf Lemma 3.2.} {\it
  Let $(M, g)$ be a hermitian manifold of complex
 dimension $m$.
Let $ \eta $ be  a $m$-form on $M$. If we write 
$ \eta = \sum_{r \geq 
0} L^{r} \eta_{r} $ with $\eta_{r}$ primitive $(m-2r)$-form,
 we have that:       
\[ \ast \eta = \sum_{r \geq 0} (-1)^{\frac{m^{2}+m}{2}+r} L^{r} C \eta_{r}, 
\;\;\;\;\;\;\;\;\;\;\;
  C \eta = \sum_{r \geq 0} L^{r} C \eta_{r}. \]                    
Then we observe that we can decompose the space of the m-forms in two parts:
                   
$\{ \eta = \sum_{r \geq 0,\, r \, even}
L^{r} \eta_{r} $ with $\eta_{r}$ primitve $(m-2r)$-form$  \}$ and 
                        
$\{ \eta = \sum_{r \geq 0, \,r \,odd}
L^{r} \eta_{r} $ with $\eta_{r}$ primitve $(m-2r)$-form$ \}$;      
                 
on the first part $ \ast =  (-1)^{\frac{m^{2}+m}{2}} C$, 
on the second one  $ \ast =  (-1)^{\frac{m^{2}+m}{2}+1} C$. }
                             
{\sc Proof.}
In \cite{Ch} p. 26 or \cite{Wel} the following formula is proved:
{\it if $\omega $ is a primitive $p$-form on $M$  and
 $r \leq m-p$ ($m= \dim _\C M$), then} 
\(\;\; 
 \ast L^{r} \omega = (-1)^{\frac{p(p+1)}{2}} \frac{r!}{(m-p-r)!} L^{m-p-r}
 C \omega . \)
      
Applying it, with $p=m-2r$, to each term $L^{r} \eta_{r} $ 
of the sum $ \eta = \sum_{r \geq 0} L^{r} \eta_{r} $,  we obtain:
\[ \ast \eta = \ast \sum_{r \geq 0} L^{r} \eta_{r} = \sum_{r \geq 
0} \ast L^{r} \eta_{r} =  \sum_{r \geq 0} (-1)^{\frac{m^{2}+m}{2}+r} L^{r}
 C \eta_{r}. \]
            
Since $C L=L C$, we have:

\vspace{0.1cm}

\hspace{5cm}
$ C \eta = \sum_{r \geq 0} C L^{r}  \eta_{r}
=\sum_{r \geq 0} L^{r} C \eta_{r}.$    
 \hfill $\Box$

\vspace{0.2cm}
                     
{\bf Corollary 3.3.} {\it
Let $(M, g)$ be a compact  K\"{a}hler  manifold
of complex dimension $m$. Consider the operators
$ \ast $ and $C$ on $H^{m}(M,\C)$.  
We have that $\ast =C$ on $K$ and $ \ast =-C$ on $K'$
(see Not. 1.2).}

\vspace{0.2cm}

Let $(M,g)$ be a compact K\"{a}hler  manifold of complex dimension $m$.
We recall the definitions of Weil's and Griffiths' Jacobians
(see for instance \cite{Ch}, \cite{Wei}, \cite{G-H}, \cite{Gre}).
                                             
$ \bullet $  Suppose that the 1-1 real 
form $\Omega $ associated to the metric is 
rational; let $p \in \N$ with $p \leq m-1$; 
               
 $p^{th}$ {\it Weil's Jacobian}   is an abelian variety so defined:
                        
the torus $H^{2p+1}(M,\R)/(H^{2p+1}(M,\Z)/torsion)$ with 
the complex structure given by $C$  and the polarization whose real part is 
${\cal R}(\alpha, \beta)   =\int \alpha \wedge \ast \beta$.

$ \bullet $ Let $p \in \N $, $p \leq m-1$;
$p^{th}$ {\it Griffiths' Jacobian} is the following complex torus:

the torus $H^{2p+1}(M,\C)/ (F^{p+1}(M) +H^{2p+1}(M,\Z))$
$=\overline{F^{n-p+1}(M)}/ \pi (H^{2p+1}(M,\Z)) $
 with the complex structure given by $i$.

\vspace{0.15cm}
             
{\sc Proof of Theorem B.}
 We recall that $H^{m}(M,\C)= K \oplus K'$  and that,              
by Corollary 3.3, $ \ast= - C$ on $K$ and  $\ast =  C$ on $K'$.
To prove the Theorem  we have only to observe that
 $H^{m}(M,\R)=(H^{m}(M,\R) \cap K) \oplus
(H^{m}(M,\R) \cap K')$
and $J_{m}=(J_{m} \cap K) \oplus  (J_{m} \cap K')$. \hfill $\Box$ 

\vspace{0.15cm}

Thus $k^{th}$-Griffiths', $k^{th}$-Weil's,
 Lazzeri's Jacobians of a K\"{a}her manifold $(M,g)$
of dimension $m= 2k+1$
 are the same real tori with a different complex structure and the 
``change'' of the complex structure depends on the complex structure of $M$
for Griffiths' J.- Weil's J. and on the class of the $(1,1)$-form associated 
to $g$ for Lazzeri's J.- Weil's J..

\vspace{0.2cm}

{\scriptsize
\begin{tabular}{|c||c|c|}
\hline 
              &  $H^{m}(M,\R)/(H^{m}(M,\Z)/torsion)$ &
 $J_{m}(M)/ \pi(H^{m}(M,\Z))$ \\ \hline \hline
$k^{th}$ G.J. &  $C|_{\{\eta + \overline{\eta }\mid\;
\eta \in J_{m}(M) \cap \overline{F^{k+1}(M)}\}} 
\oplus -C|_{\{\eta + \overline{\eta } \mid \;
\eta \in J_{m}(M) \cap F^{k+1}(M)\}}$ & 
$i|_{J_{m}(M) \cap   \overline{F^{k+1}(M)}}
 \oplus -i|_{J_{m}(M)\cap   F^{k+1}(M)}$
\\ \hline
$k^{th}$ W. J. & $C$ &  $i$  
\\ \hline 
L. J.   & $C|_{H^{m}(M,\R) \cap K} \oplus 
-C|_{H^{m}(M,\R) \cap K'}$ &  $i|_{J_{m}(M) \cap K}
\oplus -i|_{J_{m}(M) \cap K'}$
\\ \hline 
\end{tabular}
}

\vspace{0.2cm}

\vspace{0.2cm}

{\bf Definition 3.4.} {\it Let $(M, g)$ be a compact K\"{a}hler
 manfold of complex dimension $ m =2k+1$. Set 
$J'_{m}(M) := (K(M,g) \cap J_{m}(M)) \oplus \overline{(K'(M,g) \cap  
J_{m}(M))}$}.

{\bf Corollary 3.5.} {\it Let $(M, g)$ be a compact K\"{a}hler
 manfold of complex dimension $ m =2k+1$; we have 
  $T(M,g) = J'_{m}(M)/ \pi_{J'_{m}(M), \overline{J'_{m}(M)} }
 (H^{m}(M,\Z))$ with the complex structure given by 
$i$ and the imaginary part of the polarization 
$(\alpha, \beta) = - \int_{M} (\alpha + \overline{\alpha}) \wedge 
(\beta + \overline{\beta })$.}

{\bf Remark 3.6.}    
Let  $(M,g)$ be  a compact K\"{a}hler manifold of complex
dimension $m=2k+1$; seeing $T(M,g)$ as $J'_{m}(M)/ \pi(H^{m}(M,\Z))$,
 we can define a Abel map for  $k$-cycles:

let $Z_{0}$ be a $k$-cycle in $M$ and set 
{\em $B(Z_{0}) := \{ Z$   $k$-cycle  homologous to  $Z_{0}\}$};
let  $\{\psi_{1},....\psi_{l}\}$ be a basis of 
$J_{m}'(M)$:
Abel's map $\mu : B(Z_{0}) \longrightarrow T(M) $ is so defined:
if $Z \in B(Z_{0})$ and $C$ is a $(2k+1)$-chain with $\partial  C =Z-Z_{0}$,
\[\mu(Z):=(\int_{C}\psi_{1},....,\int_{C}\psi_{l}).\]

One can  easily see that the definition is good (in an analogous way as in 
\cite{Li} p.131) and                
by the same calculation   in \cite{Gri} p.826, one can  see that if
 $\{Z_{\lambda}\}_{\lambda \in B}$ is a family of effective 
$k$-subvarieties  with $B$ not singular and $Z_{0}=Z_{\lambda_{0}}$,
the map $\mu(\lambda)=\mu(Z_{\lambda})$ is not  holomorphic in general:  
let $C_{\lambda}$ be  a chain such that $\partial  C_{\lambda} 
=Z_{\lambda}-Z_{0}$; we have that $\int_{C_{\lambda}}\psi_{i}$ 
is $0$ unless $\psi_{i} $ is of kind $(k+1,k)$ or $(k,k+1)$ and that
$\int_{C_{\lambda}}\psi_{i}$ is holomorphic in the first case 
and antiholomorphic in the second case.

\vspace{0.2cm}      

{\bf Notation 3.7.}
{\it 
Let $M$ be a complex compact manifold 
 of  dimension $m$; consider a smooth deformation  of
the complex structure $ {\cal M } \rightarrow  \Delta $,
  ($\Delta $ polycylinder $\ni 0$);
we  call $ M_{t}$ the fibre over $t$ and let  $ \phi $
be a $C^{\infty}$ trivialization: ${\cal M} \rightarrow M \times \Delta$
(possibly restricting  $\Delta $); 
$\phi$ induces diffeomorphism $ \phi_{t}:$ $M_{t}\longrightarrow M$.
     Let $ \rho: T_{0}(\Delta ) \longrightarrow H^{1}( \Theta)$ be the
 Kodaira-Spencer map, where $T_{0}(\Delta )$ is the holomorphic tangent
 space to $\Delta $ in $0$ and $\Theta= \theta (T^{10}(M))$.
Suppose that $M$ is of  K\"{a}hlerian type
 (and then also $M_{t}$  by  Theorem 15 in \cite{K-S}).}
                    
We recall the definition of  Griffiths' and Weil's period maps.
        
  $\bullet $   The Griffiths period map 
${\cal G}_{q}:  \Delta \longrightarrow G_{\C}(f_{q}^{q},H^{q}(M,\C)) \times
... \times G_{\C}(f^{q-\nu }_{q} ,H^{q}(M,\C)) $
(possibly restricting $\Delta $) is the map:
 \[t \mapsto ({\cal G}^{q}_{q}(t),...,{\cal G}^{q-\nu}_{q}(t)), \]
          where 
${\cal G}^{p}_{q}(t)=  (\phi_{t}^{-1})^{\ast} F^{p}_{q}(M_{t})$,
$f_{q}^{p}= \dim F_{q}^{p}(M)$ and $\nu=[\frac{q-1}{2}]$.

 $\bullet $ If $q$ is an odd positive integer, the Weil period map
${\cal W}_{q}:  \Delta \longrightarrow G_{\C}(b_{q}/2,H^{q}(M,\C)) $
(possibly restricting $\Delta $) is the map
\[ t \, \mapsto \,  {\cal W}_{q}(t)= (\phi_{t}^{-1})^{\ast} J_{q}(M_{t}).
 \]               
Let $\phi (t) $ be a smoothly varying harmonic $(q-r,r)$ form on $M_{t}$.
Call $\phi = \phi (0)$. We have  (see \cite{Gri} p. 812 or \cite{Gre}
p. 33):
              
\vspace{0.1cm}

$  \frac{\partial \phi(t) }{\partial t}|_{t=0}
=  \rho ( \frac{\partial  }{\partial t}) \cdot \phi$,
which  is  of type  $(q-r-1,r+1)$,   \hspace{2.5cm} (1)
        
\vspace{0.1cm}
         
$  \frac{\partial \phi(t) }{\partial \overline{t}}|_{t=0}
= \overline{ \rho ( \frac{\partial  }{\partial t})} \cdot \phi$,
which is  of type $(q-r+1,r-1)$,   \hspace{2.5cm} (2)
    
\vspace{0.1cm}

where $\cdot$ is the contraction.

Thus we have that, while ${\cal G}_{q}$ is holomorphic and 
$Im \, ( d {\cal G}_{q}) (0)$ is in  $\oplus_{r=0...\nu}
Hom(H^{q-r,r},H^{q-r-1,r+1})$,
the  map ${\cal W}_{q}$ is not holomorphic; precisely, if $ \phi \in 
F^{q-r}_{q}(M)$

\vspace{0.1cm}

$ (\frac{\partial}{\partial t}{\cal G}_{q}^{q-r})(0)(\phi)=
 \pi_{{\footnotesize \overline{F_{q}^{r+1}(M)} ,F_{q}^{q-r}(M)}}
(\rho ( \frac{\partial  }{\partial t}) \cdot \phi)
    = \rho ( \frac{\partial  }{\partial t}) \cdot \phi, $
 
\vspace{0.1cm}

$ (\frac{\partial  }{\partial \overline{t}}{\cal G}_{q}^{q-r})(0)(\phi)=
\pi_{\overline{F_{q}^{r+1}(M)} ,F_{q}^{q-r}(M)}
(\overline{ \rho ( \frac{\partial  }{\partial t})} \cdot \phi)
 = 0 $;
     
\vspace{0.1cm}

while, if q is odd and  $ \phi \in J_{q}(M) $

\vspace{0.1cm}

$ (\frac{\partial  }{\partial t}{\cal W}_{q})(0)(\phi)=
 \pi_{\overline{J_{q}(M)} ,J_{q}(M)}
(\rho ( \frac{\partial  }{\partial t}) \cdot \phi)
= \rho ( \frac{\partial  }{\partial t})\cdot \phi,$ 

\vspace{0.1cm}

$ (\frac{\partial}{\partial \overline{t}} {\cal W}_{q})(0) (\phi)=
 \pi_{\overline{J_{q}(M)} ,J_{q}(M)}
  (\overline{ \rho ( \frac{\partial  }{\partial t})} \cdot \phi)
  = \overline{ \rho ( \frac{\partial  }{\partial t})} \cdot \phi $.

\vspace{0.2cm}

{\bf Corollary C.} 
{\it  Let  $(M,g)$ be  a compact K\"{a}hler  manifold   of odd complex
dimension $m$; let $\Omega \in H^{2}(M,\R)$ be the class of the $(1,1)$-form
 associated to $g$;
 consider a smooth deformation  of the complex structure
$ {\cal M } \rightarrow  \Delta $ 
for which we use Notation 3.7.
Let $\Delta' = \{ t \in \Delta \mid 
\phi_{t}^{\ast}(\Omega) \in H^{2} (M_{t}, \R)$ {\it can be represented 
by  a K\"{a}hler form for $M_{t} \}$}; suppose $\Delta'$ be a
neighbourhood of $0$.

{\bf a)} For $t \in \Delta' $ let
$g_{t}$ be  a K\"{a}hler metric on $M_{t}$ whose (1,1)-form is of class 
$\phi_{t}^{\ast}(\Omega)$; 
  if we fix a symplectic basis of $H^{m}(M,\Z)$,
we have that   the map  
$  t  \mapsto  Z(T(M_{t}, g_{t})) $
defined from  a neighbourhood of $0$ in 
$ \Delta'$ to ${\cal H}$, associating to $t$ 
 the matrix in ${\cal H}$ representing 
 $T(M_{t},g_{t})$ with that choice of the symplectic basis,
 is
holomorphic in $0$ iff $ \rho (T_{0}(\Delta)) \cdot 
(
K(M,g) \cap J_{m}(M) 
\oplus 
\overline{ (      K'(M,g) \cap J_{m}(M)    )    } )=0$
(where $ \cdot $ is the contraction), 
 (observe that $  t  \mapsto  Z(T(M_{t}, g_{t})) $
is the composition of the map $t \mapsto (\phi_{t}^{-1})^{\ast}
g_{t} $ with $\hat{T}$).

{\bf b)} Suppose now that the canonical bundle of $ M$  is trivial
and $ \rho$ is injective.
For $t \in \Delta' $ let
$g_{t}$ be  a K\"{a}hler metric on $M_{t}$ whose (1,1)-form is of class 
$\phi_{t}^{\ast}(\Omega)$, for instance
the  K\"{a}hler-Einstein metric whose (1,1)-form is of class 
$\phi_{t}^{\ast}(\Omega)$ 
($\exists !$  by  Calabi-Yau's theorem
 (see \cite{S.P.}); fix a symplectic basis of $H^{m}(M,\Z)$;
 then the differential in $0$ of    the map  
$  t  \mapsto  Z(T(M_{t}, g_{t})) $ (defined above) is injective.}
                
{\sc Proof.} We set 
${\cal W}'_{m}(t)= (\phi_{t}^{-1})^{\ast} J'_{m}(M_{t}).$     
 By (1) and (2), we have that, if $ \phi \in J'_{m}(M)$, then 

$ (\frac{\partial  }{\partial t}{\cal W}'_{m})(0)(\phi)=
 \pi_{\overline{J'_{m}(M)} ,J'_{m}(M)}
(\rho ( \frac{\partial  }{\partial t}) \cdot \phi)
= \rho ( \frac{\partial  }{\partial t})\cdot \phi,$ 

\vspace{0.05cm}

$ (\frac{\partial}{\partial \overline{t}} {\cal W}'_{m})(0) (\phi)=
 \pi_{\overline{J'_{m}(M)} ,J'_{m}(M)}
  (\overline{ \rho ( \frac{\partial  }{\partial t})} \cdot \phi)
  = \overline{ \rho ( \frac{\partial  }{\partial t})} \cdot \phi $.

See  $T (M_{t},g_{t})=  {\cal W}_{m}'(t)/ 
\pi_{  {\cal W}_{m}'(t), \overline{   {\cal W}_{m}'(t)  } } 
(H^{m}(M,\Z)) $  with  the polarization 
whose imaginary part is $(\alpha, \beta) = 
- \int_{M} (\alpha + \overline{\alpha}) \wedge (\beta + \overline{\beta })$.
         Choosen a symplectic (for $-\int \cdot \wedge \cdot$) 
basis $\{\gamma_{i}\}$ of
$H^{m}(M,\Z)$ and  a  basis $\{\omega_{j}(t)\}$
of ${\cal W}_{m}'(t)$, let
{\footnotesize $ \left( \begin{array}{c}
                                               E(t) \\
                                               F(t)      
\end{array} \right) $} be the matrix expressing $ \omega_{j} (t) $ 
in function of $ \gamma_{i}$. 
Then  $Z(T (M_{t},g_{t}))= - \overline{E(t)F(t)^{-1}}$.

 Thus $Z(T (M_{t},g_{t}))$ is  holomorphic in $t$
iff $\overline{{\cal W}_{m}'(t)} \in G_{\C}(b_{m}/2,H^{m}$ $(M, \C))$
 is holomorphic in $t$. 

If  the canonical bundle of $M$ is trivial,  the maps
$ H^{1}(\Theta) \rightarrow Hom ( H^{n,0}, H^{n-1,1})$ and
$ \overline{H^{1}(\Theta)} \rightarrow Hom ( H^{0,n}, H^{1,n-1})$ given
 by the contraction are injective (see \cite{Gri} p. 844). 
Then, as in \cite{Gri}, 
if also   $\rho$  is injective,   the maps   $( d{\cal W}_{m})(0)$ and 
$( d{\cal W}_{m}')(0)$ are injective   and we obtain b).             
\hfill $ \Box$

\section*{4. Lazzeri's Jacobian of a bundle}

{\bf Definition 4.1.} {\it
  An element $\lambda$ of a lattice $\Lambda$ is 
said primitive if  $ \not \exists \lambda' \in \Lambda$ such that 
$\lambda \in \Z\lambda'$.}

\vspace{0.1cm}

{\bf Proposition D.} {\it  
Let $(M,g_{M})$ and $(N,g_{N})$ be riemannian compact oriented manifolds
of  dimension $2(2k+1)$, resp. $2(2s)$. 
Let $p:F \rightarrow M$ be a bundle with fibre $N$
and   structure  group  $G \subset$Diff$(N)$ and suppose 
$g_{N} $ to be $G$-invariant.
We consider on $F$ the metric induced by $g_{M}$ and  $g_{N}$.   
Suppose  $\exists \;\; \lambda \in H^{2s}(N,\Z)/torsion$ such that 
$\lambda \neq 0$, $\ast_{N} \lambda = \lambda$ 
and $\lambda$ $G$-invariant. 
                                                        
We  define a map 
$e_{\lambda}:T(M) \rightarrow T(F)$
(we omit the metrics).
If $ F = M \times N$ we can define $e_{\lambda} $ simply as
the  map  induced by the map $H^{2k+1}(M,\R) \rightarrow 
H^{2k+1+2s}(M \times N,\R)$ defined by
$ \eta \mapsto \eta \wedge \lambda  $. 
More generally
define $\Lambda \in H^{2s}(F,\R)$ in the following way:
if $U$ a trivializing open subset of $M$,
$\Lambda|_{p^{-1}(U)}:= \pi^\ast \lambda$, where 
$\pi $ is the composition  of a $C^{\infty}$ trivialization
$ p^{-1}(U) \rightarrow U \times N$  with the projection
$U \times N \rightarrow N$;
 suppose $\Lambda $ is an integral form;                          
define $E_{\lambda}: H^{2k+1}(M,\R)
\rightarrow H^{2k+1+2s}(F,\R)$ 
by $E_{\lambda}(\eta)= p^{\ast} \eta \wedge \Lambda $;
  the map $E_{\lambda}$ defines a map  $e_{\lambda}:T(M) \rightarrow T(F)$.
            
The  map $e_{\lambda}$ is holomorphic and, if $ \theta_{F}$ and $ \theta_{M}$
are the polarizations of $T(F)$ and of $T(M)$, then $ e_{\lambda}^{\ast}
\theta_{F}= (\int_{N} \lambda \wedge \ast \lambda) \theta_{M}$. 
Besides $e_{\lambda}$ is  injective if one of the following conditions holds:
a) $\int_{N} \lambda \wedge \ast \lambda= 1$
b) $ F =M \times N$  and $\lambda$ is primitive
c) $ H^{\ast}(N , \Z) $ is free and $G$-invariant.}

\vspace{0.1cm}
                     

\vspace{0.1cm}
         
{\sc Proof.}                             
$\bullet $ The fact that $\ast_{N} \lambda = \lambda$ implies  at once that
    $e_{\lambda}$ is holomorphic and
 that  $e_{\lambda}^{\ast } \theta_{F} = 
(\int_{ N } \lambda  \wedge \ast \lambda) \theta_{M}$:

let $\{U_{\alpha} \}_{\alpha}$  be a trivializing covering of $M$ and let 
$\psi_{\alpha}$ be a partition of the unity for this covering; let  
$\phi_{\alpha}$ be the partition of the unity for the covering of $F$ $\{p^{-1}
(U_{\alpha})\}$  defined by $\phi_{\alpha}(y)=\psi_{\alpha}(p(y))$;     
let $\omega_{1}, \omega_{2}  \in H^{2k+1}(M,\R)$:
\hspace{-1cm}\begin{eqnarray*}
 & & \int_{F} E_{\lambda}(\omega_{1}) \wedge  E_{\lambda}(\omega_{2}) = 
 \sum_{\alpha} \int_{p^{-1}(U_{\alpha})} \phi_{\alpha} E_{\lambda}
 (\omega_{1}) 
  \wedge  E_{\lambda}(\omega_{2}) = \\
& & \sum_{\alpha} \int_{U_{\alpha} \times N} \phi_{\alpha} \omega_{1} \wedge 
  \lambda  \wedge  \omega_{2} \wedge \lambda = 
 \sum_{\alpha} \int_{U_{\alpha} } \psi_{\alpha} \omega_{1} 
 \wedge  \omega_{2} \int_{ N } \lambda  \wedge \ast \lambda = 
(\int_{ N } \lambda  \wedge \ast \lambda)
\int_{M } \omega_{1} \wedge  \omega_{2} 
\end{eqnarray*}

$\bullet$ If a) holds then the map $e_{\lambda}$ is injective
since it is  a homomorphism  of p.p.a.v.'s and preserves the polarization.
It is easy to verify that  b) implies 
 $e_{\lambda}$  injective. Finally, if c) holds, then  
 $e_{\lambda}$ is  injective by the theorem of Leray-Hirsch 
(see \cite{Sp} p.258). \hfill $ \Box$

\vspace{0.15cm}

{\bf Remark 4.2.} Suppose  $F = M \times N$.
Observe that $T(M \times N)$ is the product of the
abelian subvarieties $A_{t}$ corresponding to   
$(H^{t}(M,\R) \otimes  H^{2k+1+2s -t}(N,\R)) \oplus 
(H^{2(2k+1)-t}(M,\R)$  $ \otimes  H^{2s-2k-1 +t}(N,\R))$  
 for  $t \neq 2k+1$,  
 and  the abelian subvariety $A_{2k+1}$  corresponding to 
$H^{2k+1}(M,\R)$ $  \otimes  H^{2s}(N,\R)$.
           Obviuosly  $ Im \; e_{\lambda} \subset A_{2k+1}$.

The abelian subvarieties $A_{t}$
for $t \neq 2k+1$ are represented in the Siegel upper half
space by a matrix $Z$ with real part equal to $0$; in fact:
         
 take as a complex basis of
$(H^{t}(M,\R) \otimes  H^{2k+1+2s -t}(N,\R)) \oplus 
(H^{2(2k+1)-t}(M,\R) \otimes  H^{2s-2k-1 +t}(N,\R))$ 
a basis $\{r_{i}\} $ of  $H^{t}(M,\Z) \otimes  H^{2k+ 1 + 2s -t}(N,\Z )$;
we can complete $\{r_{i}\} $ to a  symplectic
basis $\{r_{i},s_{i}\}$ of the lattice taking  as $\{s_{i}\}$  a suitable 
 basis of $H^{2(2k+1)-t}(M,\Z) \otimes  H^{2s-2k-1 +t}(N,\Z)$;
since each $s_{i}$ is a real linear combination 
of the $\ast r_{i}'s$, $Z$ is imaginary.

\vspace{0.15cm}

{\bf Remark 4.3.} If $N$ is K\"{a}hler and the 1-1 form $\Omega$ associated to
the metric is integral, i.e. $N$ is projective, 
the hypothesis ``$ \lambda \in H^{2s}(N,\Z)/torsion$,
$\ast_{N} \lambda = \lambda$'' is satisfied if $\lambda = [\Omega^{s}]$.

 \section*{{\bf\large Appendix}}

Let  $M= \R^{2(2k+1)}/ \Lambda $ be a torus and $g_{0}$ a singular flat
 metric.             
We prove now that the definition of $T(M,g_{0})$ we gave at the end of
Introduction does not depend on the metric $g$ extending $g_{0}$:

let $\{e_{1},..., e_{n}\} $ be a basis  of ${\bf R}^{n}$ with $n =2(2k+1)$,
orthonormal for $g$ with $ < e_{1},..., e_{r} >  = (\ker g_{0})^{\perp_{g}} $  
and $ < e_{r+1},..., e_{n} > = \ker g_{0}$;

let $g'$ be another metric extending $ g_{0}$ in the  above sense;    
let $\{v_{1},..., v_{n}\} $ be a basis  of ${\bf R}^{n}$ 
orthonormal for $g'$ with $ < v_{1},..., v_{r} >  =
 (\ker g_{0})^{\perp_{g'}} $   and $ < v_{r+1},..., v_{n} > = \ker g_{0}$;     
                            
let 
$ \left( \begin{array}{cc}
                                                  A & 0 \\
                                                  C & E         
\end{array} \right) $ be        
the matrix  expressing the basis $ \{v_{i}\}$ in function of 
the basis $\{ e_{i} \}$;
                
observe that the matrix $A$ is orthogonal, in fact $ g_{0}= 
   \left( \begin{array}{cc}
                                                  I & 0 \\
                                                  0 & 0         
\end{array} \right) $ both in the basis $\{e_{ 1},...,e_{r} ,e_{r+1},
...,  e_{n}\} $ 
 and in the basis $\{v_{1},...,v_{r} , v_{r+1},...,  v_{n}\} $;
                                      
thus $ \left( \begin{array}{cc}
                                                  A & 0 \\
                                                  C & E         
\end{array} \right)^{t}   \left( \begin{array}{cc}
                                                  I & 0 \\
                                                  0 & 0         
\end{array} \right)  \left( \begin{array}{cc}
                                                  A & 0 \\
                                                  C & E         
\end{array} \right) = \left( \begin{array}{cc}
                                                  I & 0 \\
                                                  0 & 0         
\end{array} \right) $; then $A$ is orthogonal;
                             
consider the basis    $\{v'_{1},..., v'_{n}\} $ of ${\bf R}^{n}$ that  is
expressed by the matrix 
$ \left( \begin{array}{cc}
                                                  A^{-1} & 0 \\
                                                  0 & I         
\end{array} \right) $ in function of the basis 
$\{v_{1},..., v_{n}\} $;
              
the matrix expressing $\{v'_{1},..., v'_{n}\} $ in function of 
$\{e_{1},..., e_{n}\} $ is    
$ \left( \begin{array}{cc}
                                                  I  & 0 \\
                                                  C A^{-1}  & E
\end{array} \right) $;
          
thus   $\{v'_{1},..., v'_{n}\} $ is a basis orthonormal for $g'$ such that 
$<v'_{1},..., v'_{r}> = \ker g_{0}^{\perp_{g'}} $ and 
$<v'_{r+1},..., v'_{n}> = \ker g_{0} $;
              
the map  $ e_{I}^{\vee} \mapsto (v'_{I})^{\vee}$ is an isomorphism 
between $T(M, g_{0})$ built by $g$ and $T(M, g_{0})$ built by $g'$.

\vspace{1cm}
 
{\small {\bf Acknowledgements.} 
I   thank prof. F. Catanese for suggesting me to study
Lazzeri's Jacobian and for  many
helpful   discussions and suggestions on it. 
I thank also M. Manetti for some corrections to  the first part of this work.}

\end{document}